\documentclass[11pt]{amsart}

\usepackage{amsmath,amsfonts,amssymb,amscd}
\begin{document}
\renewcommand{\thefootnote}{\fnsymbol{footnote}}
\pagestyle{plain}

\title{A remark on Li-Xu's  pathology}
\author{Toshiki Mabuchi${}^*$}
\maketitle
\footnotetext{ ${}^{*}$Supported 
by JSPS Grant-in-Aid for Scientific Research (A) No. 20244005.}
\abstract
For test configurations, the Donaldson-Futaki invariant $\bar{F}_1$
 introduced in \cite{D} is well-known.
In this paper, 
a refinement $F_1$ (see~\cite{M}) of the invariant 
will be discussed.
Then we see that Li-Xu's pathology \cite{LX} doesn't occur, in the sense that 
their example of a non-normal test configuration, with vanishing $\bar{F}_1$ and trivial normalization, 
actually has non-vanishing $F_1$.
\endabstract

\section{Introduction}

In this paper, we consider a polarized algebraic manifold $(X,L)$, i.e., a pair of a nonsingular irreducible 
projective algebraic variety, defined over $\Bbb C$, and a very ample line bundle $L$ over $X$.
Put $n := \dim X$.
For the space $V = H^0(X,L)$ of holomorphic sections for $L$, we write $V$ as a direct sum 
$$
V = W \oplus W',
$$
where $W$ and $W'$ are vector subspaces of $V$ such that
the linear subsystem $W$ of the complete linear system $|L|$ is base-point free.
Then the morphism 
$$
\Phi_{W} : X \to \Bbb P^*(W)
$$ 
associated to $W$ is
easily seen to be a finite morphism onto its image $\Phi_{W}(X)$.  Assume that 
this finite morphism is not isomorphic.
For the algebraic torus $T = \Bbb C^*$,
we now consider its representation on $V$ defined by
$$
\psi : T \to \operatorname{GL}(V)
$$
such that $\psi (t)_{|W} = \operatorname{id}_{W}$ 
and $\psi (t) (W') = W'$ for all $t\in T$, while the restriction
$\psi_{|W'}$ is assumed to have only positive weights. For the complex affine line $\Bbb A^1 = \{ z\in \Bbb C\}$,
we consider the natural $T$-action on $\Bbb A^1$ by multiplication of complex numbers. 
For the Kodaira embedding 
$$
\Phi_{|L|} : X \to \Bbb P^*(V)
$$ 
associated to the complete linear system $|L|$,
we consider its image $\hat{X} := \Phi_{|L|}(X)$,
where the group $\operatorname{GL}(V)$ acts naturally on the set $\Bbb P^*(V)$ of all hyperplanes in $V$ passing through 
the origin.
Let $\mathcal{X}$ be the $T$-invariant variety
obtained as the 
closure of
$$
\bigcup_{z\in \Bbb C^*}\;\{z\} \times \psi (z) (X)
$$
in $Q := \Bbb A^1 \times \Bbb P^*(V)$, where $T$ acts on $Q$ by
$T \times Q\owns (t, (z,p)) \mapsto (tz, \psi (t) z)$.
Let $\mathcal{L}$ denote the restriction to $\mathcal{X}$ of the pullback 
$\operatorname{pr}_ 2^*\mathcal{O}_{\Bbb P^*(V)}(1)$ of the hyperplane bundle 
$\mathcal{O}_{\Bbb P^*(V)}(1)$ of $\Bbb P^*(V)$.
Here the scheme-theoretic 
fiber $\mathcal{X}_0$ of $\mathcal{X}$ over the origin in $\Bbb A^1$ 
coincides with $\{0\}\times \Phi_{W}(X)$ set-theoretically, i.e.,
$$
(\mathcal{X}_0)_{\operatorname{red}} \; =\; \{0\}\times \Phi_{W}(X),
\leqno{(1.1)}
$$
where we have $\Phi_{W}(X) \subset \Bbb P^*(W)\subset \Bbb P^*(V)$ via the projection 
of $V$ to $W$.
Then $\mu = (\mathcal{X}, \mathcal{L})$ is a typical  example of a non-normal test configuration
in Donaldson's sense \cite{D} such that the Donaldson-Futaki invariant $\bar{F}_1 (\mu )$ vanishes and that the normalization of the test configuration $(\mathcal{X},\mathcal{L})$ 
is trivial. 

\medskip
In a separate paper \cite{M}, by
considering sequences of test configurations,
we gave a refinement $F_1$ of the Donaldson-Futaki invariant $\bar{F}_1$. 
In this introduction, by a slightly different definition which is essentially equivalent 
to the one in \cite{M}, we quickly explain how to define $F_1 (\mu )$. 
For the homogeneous ideal  $I = \oplus_{\ell}\, I_{\ell}$ for $\mathcal{X}_0$ in the projective space $\Bbb P^*(V)$, we put
$$
V_{\ell} \, :=\, S^{\ell} (V)/ I_{\ell}, \qquad\ell = 1,2,\dots,
$$
where $S^{\ell} (V)$ denotes the $\ell$-th symmetric tensor product of $V$. Let $q_{\ell}$ be the Chow weight 
for $\mathcal{X}_0$ in the projective space $\Bbb P^*(V_{\ell})$.  Then 
$$
q_{\ell} \; =\;   (n+1)! c_1(L)^n[X]\,\{\,\bar{F}_1 (\mu) \ell^n + \bar{F}_2 (\mu )\ell^{n-1}+ \bar{F}_3 (\mu ) \ell ^{n-2} + \dots \,\}
$$
for $\ell \gg 1$ (see \cite{M0}). Let 
$\psi_{\ell}: T \to \operatorname{GL}(V_{\ell})$ be the representation induced by 
$\psi : T \to \operatorname{GL}(V)$. Put $N_{\ell} := \dim V_{\ell}$.
Let $\psi_{\ell}^{\operatorname{SL}}: \Bbb T_{\Bbb R}\to \operatorname{SL}(V)$ be as in (2.6).
For the weights $-\,b'_{\alpha}$, $\alpha = 1,2, \dots, N_{\ell}$,  of the $T$-action on $V_{\ell}$,
by setting $\|\psi_{\ell}\| \, :=\, \Sigma_{\alpha =1}^{N_{\ell}}\; |b'_{\alpha}|$, we define
$$
F_1 (\mu ) := 
\lim_{\ell \to \infty} \{ \,(\|\psi_{\ell} \| /\ell^{n+1})^{-1}q_{\ell}/\,\ell^n\}.
\leqno{(1.2)}
$$ 

\medskip \noindent
{\bf  Main Theorem.} 
{\em  \;\;  For $\mu = (\mathcal{X},\mathcal{L})$ above, we have $F_1 (\mu ) = -\infty$.}


\section{Proof of Main Theorem}

Let $\chi (\ell ) 
=  a_n \ell^n + a_{n-1}\ell^{n-1}+ \dots + a_1 \ell + a_0
$ be the Hilbert polynomial of 
$(X,L)$, i.e.,
$\chi (\ell ) = \dim H^0(X, L^{\otimes \ell}) $ 
for $\ell \gg 1$, where we here observe that $a_n = c_1(L)^n[X]/n!$.
We then consider the linear subspace
$$
W_{\ell} := S^{\ell}(W)/\{ I_{\ell}\cap S^{\ell}(W)\},
$$ 
of $V_{\ell}$.
For each positive real number $r$, we denote by $O(\ell^r)$ 
a function $f$ satisfying $|f| \leq C\ell^r$ for some positive real constant $C$
independent of $\ell$ and $k$, where $k$ will appear later on.
Since the linear subsystem $W$ of $|L|$ is base-point free, we obtain
$$
0\,\leq\,  N_{\ell} \,-\, n_{\ell}\, \leq \,\chi (\ell )\, -\, n_{\ell}\, =\, O(\ell^{n-1}),
\leqno{(2.1)}
$$
where $n_{\ell} := \dim W_{\ell}$.
We now observe that the weights of the $T$-representation $\psi_{\ell}$ are 
trivial when restricted to $W_{\ell}$.  Hence a basis $\{\tau_1,\tau_2, \dots, \tau_{n_{\ell}}\}$
for $W_{\ell}$ is completed to a basis $\{\tau_1,\tau_2, \dots, \tau_{N_{\ell}}\}$ 
in such a way that 
$$
\begin{cases}
&\psi_{\ell}(t)\,\tau_{\alpha} \;=\; \tau_{\alpha}\quad \;\;\;\;\text{ if $\,1\,\leq\, \alpha \leq n_{\ell}$,}\\
&\psi_{\ell}(t)\tau_{\alpha} = t_{}^{-\,b_{\alpha}}\tau_{\alpha}\,\quad \text{ if $n_{\ell}< \alpha \leq N_{\ell}$,}
\end{cases}
$$
for all $t \in \Bbb C^*$, where all $b^{}_{\alpha}$ are positive integers. Let $C_0 > 0$ be the maximum of the weights of the $T$-representation
on $V$. Then all $b^{}_{\alpha}$ satisfy $0 < b^{}_{\alpha} \leq C_0 \ell$. Hence by
(2.1), we see that
$$
0 < \gamma_{\ell}\, :=\; \Sigma_{\alpha}\, b^{}_{\alpha}/N_{\ell} \; = \;  O(1),
\leqno{(2.2)}
$$
where the summation is taken over all $\alpha$ with $n_{\ell}< \alpha \leq N_{\ell}$. Then by setting 
$$
b'_{\alpha} := \begin{cases} & \; - \gamma_{\ell},\quad\;\qquad \text{ if $\,1\,\leq\, \alpha \leq n_{\ell}$,}\\
&\;\;b_{\alpha} - \gamma_{\ell}, \quad \text{ if $n_{\ell}< \alpha \leq N_{\ell}$,}
\end{cases}
\leqno{(2.3)}
$$
we obtain $\Sigma_{\alpha =1}^{N_{\ell}}\, b'_{\alpha} = 0$. 
Moreover, by $0 < b^{}_{\alpha} \leq C_0 \ell$
and (2.2), we see from (2.3) the following:
$$
b'_{\alpha} \; = \; O(\ell ). 
\leqno{(2.4)}
$$
By considering the real Lie subgroup  $T_{\Bbb R}:= \Bbb R_+$ of $T$,
we define a Lie group homomorphism $\psi_{\ell}^{\operatorname{SL}}: T_{\Bbb R} \to \operatorname{SL}(V_{\ell})$
by
$$
\psi_{\ell}^{\operatorname{SL}}(t) \tau_{\alpha} \; =\; t^{-\,b'_{\alpha}}_{}\tau_{\alpha},
\qquad \alpha = 1,2,\dots, N_{\ell}.
\leqno{(2.5)}
$$
where $t \in \Bbb R_+$ is arbitrary. Then in view of (2.1) and (2.2), by using (2.3),
we can estimate $\|\psi_{\ell}\| := \Sigma_{\alpha =1}^{N_{\ell}}|b'_{\alpha}|$ as follows:
$$
\begin{cases}
\quad\|\psi_{\ell}\| \; &\geq \; \gamma_{\ell}\,n_{\ell} \; = \; \gamma_{\ell}\, \{\,a_n\ell^n + O(\ell^{n-1})\},\\ 
\quad\|\psi_{\ell}\| \; &\leq \; \gamma_{\ell}\,n_{\ell} + \Sigma^{N_{\ell}}_{\alpha=n_{\ell} +1}(b_{\alpha} +| \gamma_{\ell}| ) \\
& = \; 2 \gamma_{\ell} N_{\ell}\; =\; 
 2\gamma_{\ell}\, \{\,a_n\ell^n + O(\ell^{n-1})\}.
 \end{cases}
 \leqno{(2.6)}
$$
Put $\ell' := \ell k$ for positive integers $k \gg 1$.
Then the $T_{\Bbb R}$-action on $V_{\ell}$ by $\psi_{\ell}^{\operatorname{SL}}$
naturally induces the $T_{\Bbb R}$-action on $V_{\ell'}$, 
$$
T_{\Bbb R} \times V_{\ell'} \to V_{\ell'}, \qquad 
(t, \tau ) \mapsto \psi_{\ell}^{\operatorname{SL}}(t) \tau. 
\leqno{(2.7)}
$$
For a basis $\{\tau'_1,\tau'_2, \dots, \tau'_{n_{\ell'}}\}$ for $W_{\ell'}$, we complete it to a  basis
$\{\tau'_1,\tau'_2, \dots, \tau'_{N_{\ell'}}\}$ in such a way that, by the $T_{\Bbb R}$-action in (2.7), we can 
write
$$
\psi_{\ell}^{\operatorname{SL}}(t)\tau'_{\alpha}\; =\; t_{}^{-\,c_{\alpha}}\tau'_{\alpha}\,\qquad 
\alpha = 1,2,\dots,N_{\ell'},
\leqno{(2.8)}
$$
for some $c_{\alpha}\in \Bbb Q$, where $t \in \Bbb R_+$ is arbitrary. In view of (2.3) and (2.4), by comparing (2.5) 
with (2.8), we obtain
$$
\begin{cases}
&c_{\alpha} \;=\;- k \,\gamma_{\ell}\quad \;\;\;\;\text{ if $\,1\,\leq\, \alpha \leq n_{\ell'}$,}\\
&c_{\alpha}\;  =  \;  k \,O(\ell )\,\quad \,\;\;\text{ if $n_{\ell'}< \alpha \leq N_{\ell'}$.}
\end{cases}
$$
Note that by (2.1), $n_{\ell'} = a_n k^n\ell^n + k^{n-1}O(\ell^{n-1})$ and $N_{\ell'} - n_{\ell'} 
= k^{n-1}O(\ell^{n-1})$.
Hence the weight of the $T_{\Bbb R}$-action on $\wedge^{N_{\ell'}}V_{\ell'}$ 
induced by $\psi_{\ell}^{\operatorname{SL}}$ is 
$$
\Sigma_{\alpha = 1}^{N_{\ell'}} \, c_{\alpha} \; =\; -k\gamma_{\ell}n_{\ell'} + k^n O(\ell^n)\;=\;
(-\gamma_{\ell}a_n\ell^n) k^{n+1} + k^n O(\ell^n), \quad k \gg 1,
$$
so that by \cite{Mu}, Proposition 2.11, we obtain
$$
q_{\ell} \; =\; (n+1)!\cdot (-\gamma_{\ell}a_n\ell^n).
\leqno{(2.9)}
$$
By (2.6), $\|\psi_{\ell}\|/\ell^{n+1}\, \leq \, 2 \gamma_{\ell} \{\,a_n \ell^{-1} + O(1)\}$, while by (2.9), we obtain $q_{\ell}/\ell^n =  -\gamma_{\ell}a_n (n+1)!
= -\gamma_{\ell} (n+1) c_1(L)^n[X]$.
Hence we conclude from (1.2) that 
$$
F_1 (\mu ) \; =\; \lim_{\ell \to \infty} \{ \,(\|\psi_{\ell} \| /\ell^{n+1})^{-1}q_{\ell}/\,\ell^n\}
\;\leq \; \lim_{\ell \to \infty} \{-\,\ell\, (n+1)!/2\,\} \; =\; -\infty.
$$

\bigskip\noindent
{\footnotesize
{\sc Department of Mathematics}\newline
{\sc Osaka University} \newline
{\sc Toyonaka, Osaka, 560-0043}\newline
{\sc Japan}}

\end{document}